\documentclass[11pt]{article}

\usepackage[margin=1in]{geometry}
\usepackage{amsmath,amssymb, amsthm}
\usepackage{color}
\usepackage{setspace}
\usepackage{graphicx}
\usepackage[hyphens]{url}
\usepackage{subcaption}
\usepackage{multirow}
\usepackage{authblk}
\usepackage{siunitx}
\usepackage{booktabs}
\usepackage{tabularx}
\usepackage{fancyhdr}
\usepackage{float}

\newtheorem*{proposition}{Proposition}
\newcommand{\lit}{l_{i,t}}
\newcommand{\pit}{p_{i,t}}
\newcommand{\lbit}{\bar l_{i,t}} 
\newcommand{\zit}{z_{i,t}} 
\newcommand{\Unc}{\mathcal{U}} 

\title{Catastrophe Insurance: \\
An Adaptive Robust Optimization Approach}
\author[1,2]{Dimitris Bertsimas}
\author[1]{Cynthia Zeng}
\affil[1]{Operations Research Center, Massachusetts Institute of Technology}
\affil[2]{Sloan School of Management, Massachusetts Institute of Technology}
\date{}

\begin{document}

\maketitle
\setstretch{1.6}
\begin{abstract}
\noindent
The escalating frequency and severity of natural disasters, exacerbated by climate change, underscore the critical role of insurance in facilitating recovery and promoting investments in risk reduction. 
This work introduces a novel Adaptive Robust Optimization (ARO) framework tailored for the calculation of catastrophe insurance premiums, with a case study applied to the United States National Flood Insurance Program (NFIP). 
To the best of our knowledge, it is the first time an ARO approach has been applied to for disaster insurance pricing. Our methodology is designed to protect against both historical and emerging risks, the latter predicted by machine learning models, thus directly incorporating amplified risks induced by climate change. Using the US flood insurance data as a case study, optimization models demonstrate effectiveness in covering losses and produce surpluses, with a smooth balance transition through parameter fine-tuning. 
Among tested optimization models, results show ARO models with conservative parameter values achieving low number of insolvent states with the least insurance premium charged. 
Overall, optimization frameworks offer versatility and generalizability, making it adaptable to a variety of natural disaster scenarios, such as wildfires, droughts, etc. This work not only advances the field of insurance premium modeling but also serves as a vital tool for policymakers and stakeholders in building resilience to the growing risks of natural catastrophes.
\end{abstract}

\vspace{5mm}
\vspace{10mm}

\section{Introduction}

\begin{figure}
    \centering
    \includegraphics{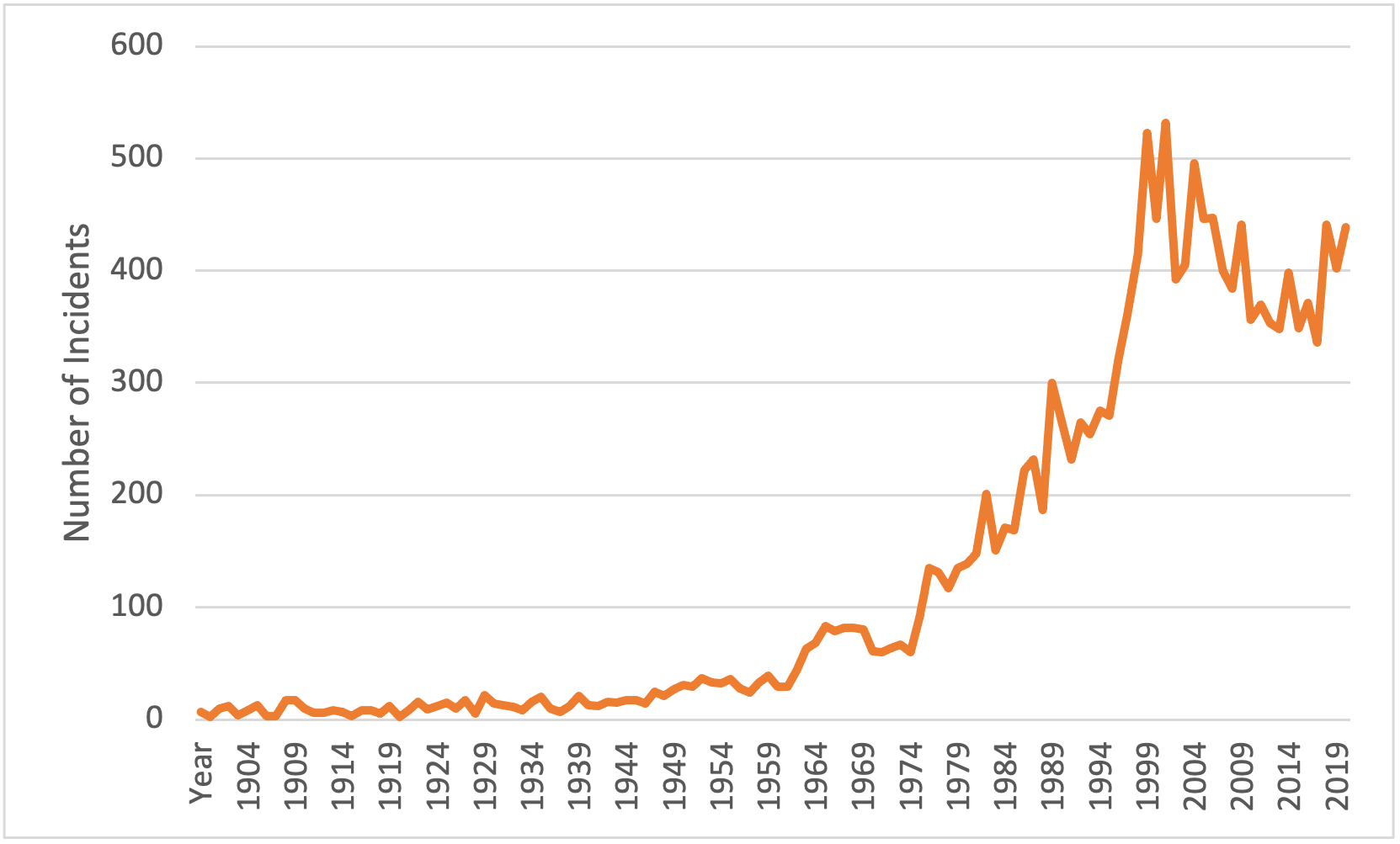}
    \caption{Number of major disasters globally since 1900, maintained by the EM-DAT database \cite{emdat}. A disaster is defined as an event which overwhelms local capacity, necessitating a request to the national or international level for external assistance. Disasters include: flood, storm, earthquake, drought, landslide, extreme temperature, wildfire, volcanic activity, mass movement (dry), glacial lake outburst, fog, etc. }
    \label{fig:disaster_trend}
\end{figure}

Global climate change causes serious consequences in climate variability and weather extremes, which could lead to more frequent and costly natural disasters worldwide \cite{van2006impacts}. As shown in Figure \ref{fig:disaster_trend}, the number of disasters worldwide has increased tremendously during the last decades. Enhanced risks for catastrophic events, such as tropical cyclones, draughts, floods, heatwaves, can inflict severe damage and losses on individuals, businesses, communities, and the entire society. It is crucial to mitigate disaster risks and facilitate climate change adaptation \cite{mercer2010disaster}. The International Panel on Climate Change (IPCC) has emphasized the need for financial instruments for disaster risk management and climate change adaptation \cite{linnerooth2015financial}. Catastrophe insurance emerges as a crucial risk management tool, offering financial support in recovery and incentivizing investments for mitigating efforts. \\

Catastrophe insurance, also known as disaster insurance, focuses on large-scale, low-frequency events that have the potential to cause widespread damage. 
 This form of insurance typically covers disasters such as hurricanes, earthquakes, floods, terrorist acts, and pandemics. The rarity of such catastrophic events complicates the insurance process, as traditional actuarial methods fall short, often due to a lack of comprehensive historical data. Compounding this challenge, climate change is leading to more regular and destructive climate-related catastrophes, which traditional reliance on historical data alone tends to underestimate.
In contrast to conventional insurance policies that spread risks across insured individuals, catastrophe insurance confronts a temporal problem of matching of regular influx of annual premiums with the irregular and unforeseeable distribution of payouts for losses.
In the United States, catastrophe insurance has historically been managed predominantly through national programs. The National Flood Insurance Program (NFIP), for example, is the principal provider for flood insurance in the country, covering more than 95\% of the underwriting risks \cite{michel2010catastrophe}. However, the NFIP's actuarial effectiveness has been the subject of scrutiny, with the program operating at a significant deficit—19 billion US dollars as of 2023—highlighting the need for reform in the structuring of such insurance schemes \cite{horn2017introduction}. The entry of private insurers into the catastrophe coverage market has been deemed crucial, yet the inherent complexities associated with rare events have led to minimal participation from private entities \cite{michel-kerjan_redesigning_2011}. Indeed, many private insurers are withdrawing from regions deemed ``uninsurable", due to escalating risks from climate change \cite{nytimes_climate_2023}. \\ 

In this work, we present an Adaptive Robust Optimization (ARO) framework for catastrophe insurance premium pricing designed to protect against uncertain losses. To the best of our knowledge, this is the first work using an ARO approach to set disaster insurance premiums. We develop the framework and implement it to flood insurance using the National Flood Insurance Program (NFIP) data. The main contributions are three-fold:
\begin{itemize}
    \item 
    
We present an Adaptive Robust Optimization (ARO) framework designed for the accurate pricing of catastrophe insurance premiums, incorporating historical data with predictions from machine learning models to offer a comprehensive view of potential losses. Specifically, we develop two types of uncertainty sets for modeling losses: a retrospective set based on historical loss distributions and a prospective set that leverages risk estimates generated by machine learning predictions.
    \item 
    We apply our ARO framework to US flood insurance using data from the National Flood Insurance Program from 1975 to 2022. Using training data, we parameterize optimization models and train our own machine learning risk models, and evaluate model performance using out-of-sample testing data. Optimization models demonstrate capabilities in effectively covering losses whilst offering the adaptability to policy-makers' risk tolerance and level of conservatism. In particular, we recommend policy makers to choose an ARO model, with conservative parameter values, to achieve superior performance in both effectiveness and efficiency in covering losses. 
    \item We highlight the adaptability and generalizability of our framework, suggesting the potential application of an ARO approach to pricing a wide range of catastrophic events, such as wildfires, droughts, extreme weather events. 
\end{itemize}

The structure of the paper is as follows. 
In Section 2, we review the relevant literature. In Section 3 we introduce the problem and outline the Robust Optimization (RO) and Adaptive Robust Optimizaiton (ARO) framework. In Section 4, we demonstrate the application of our framework through a case study on flood insurance in the United States. We explain the model parameter estimation and details on the machine learning risk prediction model. 
In Section 4, we discuss results from our models against two baselines: historical NFIP premiums and Cumulative Moving Average (CMA) scheme. Finally, we draw conclusions in Section 5.


\section{Literature Review}



This paper is related to the catastrophe modeling literature, often abbreviated at CAT modeling. It is a pivotal method for assessing and managing risks associated with extreme events employed by insurance companies. 
\cite{grossi2005catastrophe} offers a comprehensive analysis of how catastrophe models are employed for risk assessment and management purposes. Most CAT models are proprietary, with AIR Worldwide, Risk Management Solutions (RMS), and EQECAT being the major players in the private sector; and the open-source HAZUS model developed by FEMA \cite{schneider2006hazus}. 
One major challenge of the CAT modeling approach lies in the lack of historical data due to the rarity of the events, and thus standard actuarial techniques fail to capture tail-event risks.
In addition, catastrophe modeling depends on scenario simulations, which can be both numerous and lead to largely varying outcomes. How to combine different``what-if" scenarios remains a challenge to decision makers \cite{kunreuther_risk_2013}. 
Our paper illustrates the promise of an optimization-based approach to model uncertain losses and climate risks directly, thus offering transparency and offer greater robustness against rare events.\\ 


Our work is also related to the literature of disaster management. Many works have focused the problem of general resource allocations to different programs or regions under a given budget constraint \cite{wang2023risk, alem2016stochastic, salmeron2010stochastic}. \cite{yang2021floodflashfund} discusses fund allocation for flash flood reduction, and \cite{liu2022funding} further incorporates the use of insurance premiums as a source of funding. As highlighted by \cite{kunreuther_risk_2013}, the critical issue in insurance policy lies in the need for decision-making robustness in the face of climate change's uncertainties.
There is also growing work on using robust optimization in disaster relief management to deal with uncertainties \cite{ben2011robust, zokaee2016robust}. 
Few works discuss the use of optimization for catastrophe insurance pricing directly, \cite{ermolieva_flood_2017} applies stochastic optimization to the Dutch Flood insurance scheme, and
\cite{ermolieva2013modeling} discusses an integrated catastrophe management approach by incorporating CAT models with stochastic optimization methods.  
Our study fills in the gap by proposing a robust optimization framework by modeling uncertain losses and integrating machine learning forecast risks to address the increasing unpredictability of weather events driven by climate change.\\


Finally, this work broadly belongs to the climate finance literature, and see \cite{hong_climate_finance_2020} for a comprehensive overview. 
\cite{bouwer_financing_climate_2006, hochrainer-stigler_funding_2014} examine the public policy implications in funding climate change adapation. 
It is critical to design new tools to financially manage weather risks. 
\cite{mills_synergisms_2007} discusses the role of insurance sector in decreasing the vulnerability of human and natural systems.
\cite{auh_municipal_bond_2022, fowles_accounting_2009} propose the use of municipal bonds to finance natural disasters. 
\cite{keucheyan2018_insuring_climate_change} considers insuring climate change induced risks across broad business spectrum. 
Our paper adds to the literature by proposing a new design of the weather-related insurance contract to manage weather risk.

\section{Optimization Framework for Catastrophe Insurance}
In this section, we introduce the problem and present the Robust Optimization (RO) and Adaptive Robust Optimizaiton (ARO) frameworks.
\subsection{Robust Optimization}
Robust Optimization (RO) is a useful methodology for handling optimization problems with uncertain data \cite{Bertsimas2011}. It has been applied to address uncertainties in various fields, including operations research, engineering, and finance. Consider a generic linear programming problem
\[
\min _{\mathbf{x}}\left\{\mathbf{c}^T \mathbf{x} \mid \mathbf{A} \mathbf{x} \leq \mathbf{b}\right\} ,
\]
where $\mathbf{c} \in \mathbb{R}^n, \mathbf{b} \in \mathbb{R}^m$ and $\mathbf{A} \in \mathbb{R}^m \times \mathbb{R}^n$. Robust Optimization addresses the problem where data (c, $\mathbf{A}, \mathbf{b})$ are uncertain, but are known to reside in an uncertainty set $\mathcal{U}$. Based on prior information or assumptions, we construct the uncertainty set $\mathcal{U}$ to express the uncertainties in data.  We are addressing a family of problems for each realization of $(\mathbf{c}, \mathbf{A}, \mathbf{b}) \in \mathcal{U}$. Therefore, we can reformulate the problem into its Robust Counter part
\[
\min _{\mathbf{x}}\left\{\mathbf{c}^T \mathbf{x} \mid \mathbf{A x} \leq \mathbf{b} \quad \forall(\mathbf{c}, \mathbf{A}, \mathbf{b}) \in \mathcal{U}\right\}. 
\]

\subsection{Nominal Formulation} 

We consider setting insurance premiums for N locations for the insurance period of T years, which we denote with variable $\pit$, where $i=1, \dots, N$, $t = 1, \dots, T$. We are given historical losses for each location for each year in the past $T_0$ years, which we denote with $ \lbit $, where $i=1, \dots, N$, $t = 1, \dots, T_0$. We assume that the future losses for each location $i$ n the insurance period of T years with $\lit$, $i=1, \dots, N$, $t = 1, \dots, T$. Note that this quantity is unknown, but we assume the knowledge of it to introduce a simple deterministic LP formulation to introduce the basic requirements and set up the generic framework. In the next subsections, we expand on how to model this uncertain quantity through Robust Optimization framework. \\ 

We formulate an LP model to set insurance premium price. The objective function is to minimize the overall premium collected, or the least required premium needed. In addition, to model the consumer behavior that as we increase premium price less consumers are willing to purchase the insurance product, we introduce a damping function $f: \mathbb{R} \to \{0,1\}$, a monotonically decreasing function representing decline in demand due to higher premiums. Details of the choice of the damping function will be discussed later in Section 4.2.2. The overall objective is as follows 
\begin{equation}
    \min_{\pit} \sum_{i=1}^N \sum_{t=1}^T 
    f(p_{i,t}) * p_{i,t}. 
\end{equation}

We require the premium price to cover projected losses with an additional buffer amount, denoted by $\delta$, which is set at constant for each location. Thus we require 
\begin{equation}
    \sum_{t=1}^T 
    f(p_{i,t}) * p_{i,t}  - \sum_{t=1}^T f(p_{i,t}) * l_{i,t} \ge \delta,  \quad i \in [N]. 
\end{equation}
In addition, we impose a constraint to require premiums collected over consecutive years to vary slowly, in order to prevent drastic changes in insurance premiums 
\begin{equation} \label{eq:lp_coverloss}
        | p_{i,t} - p_{i, t-1} | \le \gamma_1,
    \quad i \in [N], \quad t \in [T]. 
\end{equation}
Variables $\pit$ should be positive, for each location $i$ for each period $t$ 
\begin{equation}
\pit \in \mathbb{R}^+, \quad i \in [N], \quad t \in [T]. 
\end{equation}

\subsection{Robust Optimization Formulation}
In the nominal formulation, we assume the knowledge of projected loss for each of the future period. However, this quantity is unknown and highly uncertain. In this subsection, we expand on how to construct uncertainty sets to describe this quantity. \\ 

The overall optimization formulation is the same as before, except for constraint \ref{eq:lp_coverloss}, where we require the inequality to hold for all losses in the uncertainty set 
\begin{equation}
    \sum_{t=1}^T 
    f(p_{i,t}) * p_{i,t}  -  \sum_{t=1}^T f(p_{i,t}) * l_{i,t} \ge \delta,  \quad i \in [N], \quad \forall \lit \in \mathcal{U}. 
\end{equation}

We propose two uncertainty sets to model the uncertainties of the future losses: with Central Limit Theorem (CLT) and with Machine Learning driven risks. 

\subsubsection{Uncertainty Set from Central Limit Theorem}
Based on the assumption that for each specific location, future flooding losses follow the distribution of historical losses. We adopt the central limit theorem (CLT) to form the uncertainty set as discussed in \cite{bertsim_ro}. Formally, for a fixed location $i$, $L_{i,t}, t \in [T]$ are independent, identically distributed random variables with mean $\bar{l_i}$ and standard deviation $\bar{\sigma_i}$, where $\bar{l_i}, \bar{\sigma_i}$ are historical mean and standard deviation for location $i$. We assume the uncertain quantities $L_{i,t}$ take values such that 
\begin{equation}
    \left|\sum_{t=1}^T L_{i,t}- T \cdot  \bar{l}_{i} \right| \leq   \gamma_2 \cdot \sigma_i \sqrt{n}, 
\end{equation}
where $\gamma_2$ is a small constant to denote how close future losses distribution should deviation from the normal distribution. In other words, we describe the uncertain quantities $L_{i,t}$ as values in the uncertainty set 
\begin{equation}
    \mathcal{U}^{CLT}_i = \left \{ \left(l_{i,1}, \ldots, l_{i,T} \right) : \frac{ | \sum_{t=1}^{T}l_{i,t} - \bar{l}_{i} | } {\sigma_i \sqrt{T}} \le \gamma_2 \right \} ,  
\end{equation}

where $\bar{l}_{i}, \sigma_i$ can be computed for each location $i$ using historical data. The larger we set $\gamma_2$, the more conservative the optimization model, and higher premiums will be. As a remark, we derive one uncertainty set for each location using the historical mean and standard deviation for that location to acomodate different flooding risk profiles. 
We derive the uncertainty sets separately for each location $\mathcal{U}^{CLT}_i $.

\subsubsection{Uncertainty Set from Machine Learning Risk Models}
In addition, suppose we have some information on how future losses should be, which could be informed by risk models. We can formulate such model predictions as uncertainty sets. One way of obtaining such risk prediction is through machine learning models, as recent advances in ML models demonstrate capabilities for models to predict accurate multi-year forecasts. In this work, we build machine learning models to obtain flooding risks, see greater details in the next section. Nevertheless, one can obtain such risk predictions from physical-based models, or other alternative approaches.  \\ 

We define a major flood event to be flood incurring loss over a threshold loss level $\Theta$. Suppose we have model predictions for the risk that at location $i$, the probability of having one flood event in the next $k$ years to be $q_{i,k}$. Then we can express the model prediction as 
\begin{equation}
   \mathbb{P} \left( \sum_{t=1}^{k} z_{i,t} =1 \right) = q_{i,k},  
\end{equation}
where $z_{i,t} \in \{ 0,1\}$  are binary variables denoting if there is a major flood at location $i$  at time $t$.
By modeling such, we assume having more than one major event is negligible for reasonable $k$.
Since model predictions are probabilistic predictions which can have errors, and we want to be conservative and protect against suffering potential huge losses. Thus, assuming actual incidence of having a major flood is close to the model predictions, we can express $\zit$ as random variables taking values in an uncertainty set as  
\begin{equation}
\mathcal{U} = \left \{ z_{i,t} : | \sum_{t=1}^{k} z_{i,t} - q_{i,k} |\leq \epsilon , 
\quad  \sum_{t=1}^{k} z_{i,t} \leq 1  \right \} ,
\end{equation}
where the constant $\epsilon$ is a parameter to indicate how close we believe the model predictions are to actual probabilities. The larger the value, the less confident and more conservative our model will be. Therefore, linking variables $z_{i,t}$ to $\lit$, and assuming the future losses will be bounded by the expected value coming from suffering a major flood, we model the uncertainty set of the future loss for each location $i$ as follows 

\begin{equation} \label{eq:ro_ml_u_set}
\mathcal{U}_i^{ML} = \left\{ z_{i,t}, l_{i,t} : 
    \sum_{t=1}^{k} \lit \leq \Theta \cdot \sum_{t=1}^{k} z_{i,t}, \quad 
    | \sum_{t=1}^{k} z_{i,t} - q_{i,k} | \leq \epsilon, \quad 
    \sum_{t=1}^{k} z_{i,t} \leq 1 
    \right\} .
\end{equation}

\subsubsection{The Robust Counterpart}
Combining uncertainty sets from Central Limit Theorem and from Machine Learning risk models, the robust optimization formulation of the problem is as follows 
\begin{align} \label{eq:ro_form}
\begin{split}
\min_{\pit} & \sum_{i=1}^N \sum_{t=1}^T 
    f(p_{i,t}) * p_{i,t}  \\
&   \sum_{t=1}^T 
    f(p_{i,t}) * p_{i,t}  - \sum_{t=1}^T f(p_{i,t}) * l_{i,t} \ge \delta,   \quad \forall \lit, \zit \in \mathcal{U}_i^{CLT}, \lit \in \mathcal{U}_i^{ML}, \quad i \in [N], \\
&      || p_{i,t} - p_{i, t-1} || \le \gamma_1,
    \quad i \in [N], t \in [T],  \\ 
& \pit \in \mathbb{R}^+ ,\zit \in \{0,1\}, \quad i \in [N], t \in [T] .  
\end{split}
\end{align}

Since the uncertain occurs only in one constraint, we can write Problem \ref{eq:ro_form} as a min-max problem 
\begin{align} \label{eq:ro_form2}
\begin{split}
\min_{\pit} \max_{\lit, \zit \in \Unc} & \sum_{i=1}^N \sum_{t=1}^T 
    f(p_{i,t}) * p_{i,t}  \\
&   \sum_{t=1}^T 
    f(p_{i,t}) * p_{i,t}  - \sum_{t=1}^T f(p_{i,t}) * l_{i,t} \ge \delta,  \quad i \in [N] , \\
&      || p_{i,t} - p_{i, t-1} || \le \gamma_1,
    \quad i \in [N], t \in [T] , \\ 
& \lit \in \mathcal{U}_i^{CLT},  \quad i \in [N] \\
& \lit, \zit \in \mathcal{U}_i^{ML},  \quad i \in [N], \\
& \pit \in \mathbb{R}^+, \zit \in \{0,1\} \quad i \in [N], t \in [T] .    
\end{split}
\end{align}

Next, we consider the inner problem
\begin{align} 
\begin{split}
\max_{\lit, \zit \in \Unc} & \sum_{i=1}^N \sum_{t=1}^T 
    f(p_{i,t}) * p_{i,t}  \\
&   \sum_{t=1}^T 
    f(p_{i,t}) * p_{i,t}  - \sum_{t=1}^T f(p_{i,t}) * l_{i,t} \ge \delta,  \quad i \in [N],\\
& \lit \in \mathcal{U}^{CLT}, \\
& \lit, \zit \in \mathcal{U}^{ML}, \\
& \pit \in \mathbb{R}^+, \zit \in \{0,1\} \quad i \in [N], t \in [T] .    
\end{split}
\end{align}
For the inner problem, we can treat $\pit$ as constants, thus the inenr problem can be simplified to 
\begin{align} 
\begin{split}
\max_{\lit, \zit \in \Unc} & \sum_{i=1}^N \sum_{t=1}^T l_{i,t} \\ 
& \lit \in \mathcal{U}^{CLT}, \\
& \lit, \zit \in \mathcal{U}^{ML} ,\\
& \pit \in \mathbb{R}^+, \zit \in \{0,1\} \quad i \in [N], t \in [T]  .   
\end{split}
\end{align}
Note that the uncertainty set is composed of two separate uncertainty sets $\Unc^{CLT}$ and $\Unc^{ML}$. We can thus decompose the inner problem into two subproblems, and the objective of the original problem takes the maximum of the two subproblems. Solving each subproblem separately, we can then plug back the analytical solution from each subproblem back to the original problem. 

\begin{proposition}
The overall min-max problem is equivalent to
\begin{align}\label{eq:ro_reform2}
\begin{split}
    & \min_{\pit} \sum_{i, t}^T
    f(\pit) * \pit  \\
    &\sum_{t=1}^{T} 
    f(\pit) * \pit  - \frac{1}{T}\sum_{t=1}^T f(p_{i,t}) * L^{CLT}_i  \ge \delta, \quad \forall i \in [N] ,\\
    &\sum_{t=1}^{k} 
    f(\pit) * \pit  - \frac{1}{k}\sum_{t=1}^k f(p_{i,t}) * L^{ML}_i  \ge \delta, \quad \forall i \in [N] ,\\
    &|| \pit  - p_{i, t-1} || \le \gamma_1,
    \quad \forall i \in [N], \quad t\in [T], 
\end{split}
\end{align}
where 
\begin{align} \label{eq:L_CLT}
    L_i^{CLT} &= T\cdot \bar{l}_{i} + \gamma_2 \cdot \sigma_i \sqrt{T}, \\ \label{eq:L_ML}
     L^{ML}_i &=  \Theta \cdot \min  \{1, q_{i,k}  + \epsilon \} .
\end{align}

\end{proposition}
As a remark, we take the average "damping" over losses  $\lit$ because we solve for the optimal $\sum_t \lit$. 

\begin{proof}
 Consider the first subproblem
\begin{align}
\begin{split}
     \max_{\lit} & \sum_{i=1}^N \sum_{t=1}^T l_{i,t} \\ 
    & | \sum_{t=1}^{T} \lit - T \cdot\bar{l}_{i} |   \le \gamma_2 \cdot \sigma_i \sqrt{T} \quad   \forall i \in [N],\\
    &\lit \in \mathbb{R}^+,
\end{split}
\end{align}
looking at the constraint, we can take out  $| \cdot |$ since we are maximizing over $\lit$, and rearranging terms 
\begin{align}
     \sum_{t=1}^{T} \lit \le T \cdot\bar{l}_{i}  + \gamma_2 \cdot \sigma_i \sqrt{T} \quad   \forall i \in [N],
\end{align}
and optimality is achieved at equality, we can solve for i.e., $\lit^*$ s.t. 
\begin{align}
 &\sum_{t=1}^{T} \lit  = T \cdot\bar{l}_{i}  + \gamma_2 \cdot \sigma_i \sqrt{T} \quad \forall i \in [N] .
\end{align}
Hence we can compute for all location $i \in [N] $, and denote the analytical solution to the first subproblem as $L_i^{CLT}$ 
\begin{equation} \label{eq:ro_reform_sol1}
    L_i^{CLT} := \sum_{t=1}^{T} \lit^* =  T \cdot\bar{l}_{i}  + \gamma_2 \cdot \sigma_i \sqrt{T} \sqrt{(T_1 - T_0)}. 
\end{equation}

Consider now the second subproblem 
\begin{align}
\begin{split}
   \max_{\lit \zit} & \sum_{i=1}^N \sum_{t=1}^T l_{i,t}\\
& \sum_{t=1}^{k} \lit \leq \Theta \cdot \sum_{t=1}^{k} z_{i,t} , \\ 
&| \sum_{t=1}^{k} z_{i,t} - q_{i,k} | \leq \epsilon , \\ 
&\sum_{t=1}^{k} z_{i,t} \leq 1 , \\
& \lit \in \mathbb{R}^+, \zit \in \{0,1\}   .
\end{split}
\end{align}

Without loss of generality, we relax $\zit$ to take continuous value $\zit \in [0, 1]$, because we can treat $\sum_i^{k} z_{i,t}$ as one variable taking continuous values in $[0,1]$. Thus looking at the constraints concerning $\zit$, and take out  $| \cdot |$ since we are maximizing over $\zit$, we get $\zit^*$ achieve optimality at 
\begin{equation}
    \sum_{t=1}^{k} z_{i,t}^* = \min  \{1, q_{i,k}  + \epsilon \} ,
\end{equation}

which gives $\lit^*$ at optimality at 
\begin{equation}
    \sum_{t=1}^{k} \lit^* = \Theta \cdot \sum_{t=1}^{k} z_{i,t}^* .
\end{equation}

Combining with the solution from the first subproblem given by equation \ref{eq:ro_reform_sol1} we have the sum of future loss for each location given as follows 
\begin{align}
    \sum_{t=1}^{T} \lit^*  &=  T \cdot\bar{l}_{i}  + \gamma_2 \cdot \sigma_i \sqrt{T} ,\\ 
    \sum_{t=1}^{k} \lit^*  &=  \Theta \cdot \min  \{1, q_{i,k}  + \epsilon \} ,
\end{align}
where the first equality bounds the entire future period $T$, the second equality bounds the period depending on machine learning model's risk forecast horizon k. As a remark, we can have multiple risk models for different forecasting horizons. 
\end{proof}

\subsection{Choices of demand damping function f(p)}
Recall that to model the behavioral aspect that as insurance premium increase there is less demand for it, we introduce the damping function $f(p): \mathbb{R} \to \{0,1\}$  into the constraint. 
To preserve the convexity property of the overall problem, we model the behavior using piece-wise linear functions under the scope of this work. 
In the next section, we explain in greater detail on the estimation using NFIP data, as well as sensitivity analysis on the choice of the function. \\

As a remark, as we damp demand, we correspondingly damp the covered loss in the constraint. Since the inner problem gives the the overall loss over the forecasting period $T$, i.e., $L_i^{CLT}$ gives the maximum loss deduced from the uncertainty set over period $T$, we have taken the corresponding damping term to be the average over $T$ periods, i.e., $\frac{1}{T}\sum_{t=1}^T f(p_{i,t})$.

\subsection{Adaptive Robust Optimization Formulation} 


\newcommand{\ait}{\alpha_{i,t}}
\renewcommand{\bit}{\beta_{i,t}}


This section discusses how to the RO framework with adaptive robust optimization techniques, enabling premium adjustments based on actual loss experiences. We let the premiums depend on realized losses using affine decision rules, as proposed in Chapter 7 of \cite{bertsim_ro}. This approach not only refines premium pricing accuracy but also ensures a responsive and equitable insurance mechanism against the backdrop of unpredictable catastrophic events. In particular, we let premiums depend on loss from the previous time step as follows  

\begin{equation} \label{eq:aro_rule}
p_{i,t} = 
\begin{cases}
\alpha_{i,1} ,& \text{for } t = 1 \\
\alpha_{i,t} + \beta_{i,t} \cdot l_{i,t-1} ,& \text{for } t = 2, \dots, T
\end{cases}
\end{equation}
where premium for location $i$ at time period $t$ is determined by a linear combination of parameters. Specifically, for the first time period, the premium is set to a base value $\alpha_{i,1}$; for subsequent periods, the premium is adjusted based on the loss $l_{i,t-1}$, experienced in the previous period, with $\ait$ and $\bit$ new variables to be optimized over.\\ 
 
Throughout this section, we consider the simplified robust optimization problem from the original problem, by dropping the demand-damping to allow the derivation of a robust counter part. As discussed later in the sensitivity analysis in Section \ref{Sensitivity}, the demand damping term does not materially influence the premium outcome. Thus the problem we consider is given as follows 
\begin{align} 
\begin{split}
\min_{\pit} & \sum_{i=1}^N \sum_{t=1}^T 
     p_{i,t}  \\
&   \sum_{t=1}^T 
    p_{i,t}  - \sum_{t=1}^T  l_{i,t} \ge \delta,   \quad \forall \lit \in \mathcal{U}_i^{CLT}, \quad i \in [N] ,\\
&      || p_{i,t} - p_{i, t-1} || \le \gamma_1,
    \quad i \in [N], \quad  t \in [T] , \\ 
& \pit \in \mathbb{R}^+ ,\zit \in \{0,1\}, \quad i \in [N], \quad t \in [T]. 
\end{split}
\end{align}
Recall that this optimization problem aims to minimize the total premiums over all locations and time periods while ensuring that the cumulative premium exceeds the cumulative losses by at least a margin of $\delta$. Note that in this formulation, we impose slowly varying constraint on $\ait$, instead of $\pit$ as in the previous formulation, because $\pit$ is now depending on uncertain variables $\lit$.\\ 

 Similar as before, noticing that the constraints for each location $i$ is independent of other locations, and thus minimizing the aggregated premium is the same as minimizing the premium for each location. Therefore we can decompose the problem by solving for each location $i$ the ARO problem independently. In addition, noticing that the objective function now depends on uncertain variable $\lit$, we therefore use the epigraph formulation as discussed in Chapter 2 of the book \cite{bertsim_ro} to move all variables containing uncertain variables into constraints. We arrive at the following ARO formulation for one location $i$
\begin{align} \label{eq:aro_form}
\min_{\ait, \bit} \Omega\\ 
& \sum_{t=1}^T 
     p_{i,t} \leq \Omega,   \quad \forall l_{i,t} \in \mathcal{U}_i^{CLT},
     \label{eq:aro_cons1}\\
&   \sum_{t=1}^T 
    p_{i,t}  - \sum_{t=1}^T  l_{i,t} \ge \delta,   \quad \forall l_{i,t} \in \mathcal{U}_i^{CLT} 
    \label{eq:aro_cons2} ,\\
&      || \alpha_{i,t} - \alpha_{i, t-1} || \le \gamma_3,
    \quad  t = 2, \dots ,T  ,\\ 
&     || \beta_{i,t} - \beta_{i, t-1} || \le \gamma_4,
    \quad  t = 2, \dots ,T , \\ 
& \pit \geq 0, \quad  t = 1,\dots,T, \quad \forall l_{i,t} \in \mathcal{U}_i^{CLT}.\label{eq:aro_cons3} , 
\end{align}
where $\pit$ is a quantity that depends on the uncertain losses given by equation \ref{eq:aro_rule}. Next, we derive the robust counter part by taking the RC for each constraint, that depends on uncertain variable $\lit$, independently. \\ 

First, consider constraint give by inequality \ref{eq:aro_cons1}, substituting $\pit$ with $\ait, \bit$  which now depends on the uncertain variable $\lit$. 
We rewrite the constraint as follows 
\begin{align}
\sum_{t=1}^{T} \ait + \sum_{t=2}^{T} \bit \cdot l_{i,t-1} \le \Omega ,   \quad \forall l_{i,t} \in \mathcal{U}_i^{CLT},
\end{align}
which is is equivalent to
\begin{equation}
    \sum_{t=1}^{T} \ait + \max_{\lit \in \mathcal{U}_i^{CLT}} \left\{ \sum_{t=2}^{T} \bit \cdot l_{i,t-1} \right\} \leq \Omega. 
\end{equation}

Consider now the inner maximization problem using the explicit expression for $\Unc_i^{CLT}$, we have a LP problem in $\lit$
\begin{align}
\begin{split}
\max_{\lit} & \sum_{t=2}^{T} \bit \cdot l_{i,t-1}\\ 
\sum_{t=1}^{T} \lit &\leq \bar{T}_i + \gamma_2 \cdot \sigma_i \cdot \sqrt{T} ,\\
-\sum_{t=1}^{T} \lit &\leq -\bar{T}_i + \gamma_2 \cdot \sigma_i \cdot \sqrt{T}.  
\end{split}
\end{align}
The inner problem thus satisfies strong duality. Notice that the uncertainty set is given in polyhedron form, we can thus take the dual by introducing dual variables $s^1_1, s^1_2$, and arrive at the following form 
\begin{align}
\begin{split}    
\min_{s^1_1, s^1_2} \sum_{j=1}^{2} c_j s^1_j \\ 
s^1_1 - s^1_2 &\geq \beta_t, \quad \forall t=2,\dots,T, \\
s^1_1 - s^1_2 &\geq 0, \\
s^1_1, s^1_2 &\geq 0 . 
\end{split}
\end{align}

By strong duality, the inner maximization problem has the same objective value of the dual minimization problem. Using the dual expression, constraint \ref{eq:aro_cons1} becomes 
\begin{align*}
\sum_{t=1}^{T} \ait + \min_{s^1_1, s^1_2} \sum_{j=1}^{2} c_j s_j \leq \Omega\\
s^1_1 - s^1_2 &\geq \bit, \quad \forall t=2,\dots,T ,\\
s^1_1 - s^1_2 &\geq 0 ,\\
s^1_1, s^1_2 &\geq 0. 
\end{align*}

Note that we can take away the minimization term because if any feasible $s^1_1, s^1_2$ satisfies this constraint, the minimum also does. \\ 

Second, consider constraint give by inequality \ref{eq:aro_cons2}, which ensurs premiums cover losses 
\begin{equation}
     \sum_{t=1}^T 
    p_{i,t}  - \sum_{t=1}^T  l_{i,t} \ge \delta,   \quad \forall l_{i,t} \in \mathcal{U}_i^{CLT} . 
\end{equation}
Multiplying both sides by $-1$, we have 
 \begin{equation}
 \sum_{t=1}^{T} \lit  - \sum_{t=1}^{T} p_{it} \leq -\delta. 
 \end{equation}
Substituting $\lit$ with $\ait$ and $\bit$ 
\begin{equation}
\sum_{t=1}^{T} \lit  - \left[ \sum_{t=1}^{T} \ait + \sum_{t=2}^{T} \bit l_{i,t-1} \right] \leq -\delta,
\end{equation}
re-arranging terms to collect all the uncertain terms $\lit$
\begin{equation}
    -\sum_{t=1}^{T} \ait + \left[ l_{i,1}  + \sum_{t=2}^{T} (1 - \bit)\lit  \right] \leq -\delta,
\end{equation}
we finally arrive at  
\begin{equation}
    -\sum_{t=1}^{T} \ait + \max_{\lit \in \mathcal{U}_i^{CLT}} \left\{ l_{i,1}  + \sum_{t=2}^{T} (1 - \bit)\lit  \right\}  \leq -\delta. 
\end{equation}
Consider now the inner maximization problem, which is once again a LP in $\lit$
\begin{align}
\begin{split}
\max_{\lit}
& \quad l_{i,1}  + \sum_{t=2}^{T} (1 - \bit)\cdot \lit \\ 
\sum_{t=1}^{T} \lit &\leq \bar{T}_i + \gamma_2 \cdot \sigma_i \cdot \sqrt{T} , \\
-\sum_{t=1}^{T} \lit &\leq -\bar{T}_i + \gamma_2 \cdot \sigma_i \cdot \sqrt{T}. 
\end{split}
\end{align}
Similar as before, by strong duality, the dual is given as a minimization problem in dual variables $s^2_1, s^2_2$
\begin{align}
\begin{split}    
\min_{s^2_1, s^2_2} \sum_{j=1}^{2} c_j s^2_j \\ 
s^2_1 - s^2_2 &\geq 1- \bit, \quad \forall t=2,\dots,T ,\\
s^2_1 - s^2_2 &\geq 1 ,\\
s^2_1, s^2_2 &\geq 0 . 
\end{split}
\end{align}
 Therefore constraint \ref{eq:aro_cons2}  becomes 
 \begin{align}
     \begin{split}
 & -\sum_{t=1}^{T} \ait + \sum_{j} c_j s^2_j \leq -\delta ,\\ 
& s^2_1 - s^2_2 \geq 1 - \bit, \quad \forall t=2,\dots,T ,\\ 
& s^2_1 - s^2_2 \geq 1, \\ 
& s^2_1, s^2_2 \geq 0. 
     \end{split}
 \end{align}
Finally, we consider the positivity constraint as given by inequality \ref{eq:aro_cons3}
\begin{equation}
    \pit \ge 0, \quad \forall t \in [T], \quad \forall \lit \in \Unc^{CLT}
\end{equation}
substituting $\pit$ with $\ait, \bit$, the constraint is equivalent to 
\begin{equation}
    \ait + \bit \cdot l_{i, t-1} \ge 0, \quad \forall t \in [T], \quad \forall \lit \in \Unc^{CLT}
\end{equation}
Since for each time $t$, we have one separate constraint. Therefore, considering this constraint for some t 
\begin{equation*}
\ait + \beta_{i, t} \cdot l_{i, t-1} \geq 0, \quad \forall \lit \in \Unc^{CLT}. \\ 
\end{equation*}
Multiplying both sides by $-1$ and rearranging terms, we have 
\begin{equation*}
    -\ait + \max_{\lit \in U} \left\{ - \beta_{i, t} \cdot l_{i, t-1} \right\} \leq 0 .
\end{equation*}

Similar as before, we consider the inner maximization problem 
\begin{equation}
    \max_{\lit \in U} - \beta_{i, t} \cdot l_{i, t-1} ,
\end{equation}
which has strong duality, with dual given for each t in dual variables ${ s^{3,t}_1,  s^{3,t}_2}$ as 
 \begin{align}
     \begin{split}
 & \min_{ s^{3,t}_1,  s^{3,t}_2} \sum_{j} c_j s^{3,t}_j  \\ 
& s^{3,t}_1 - s^{3,t}_2 \geq - \beta_{i, t} ,\\ 
& s^{3,t}_1 - s^{3,t}_2 \geq 0 ,\\ 
& s^{3,t}_1, s^{3,t}_2 \geq 0. 
     \end{split}
 \end{align}

 Thus the positivity constraint for some $t$ becomes 
 \begin{align}
     \begin{split}
 & \ait + \sum_{j} c_j s^{3,t}_j \leq 0 \\ 
& s^{3,t}_1 - s^{3,t}_2 \geq - \beta_{i, t} , \\ 
& s^{3,t}_1 - s^{3,t}_2 \geq 0 ,\\ 
& s^{3,t}_1, s^{3,t}_2 \geq 0. 
     \end{split}
 \end{align}

 Thus, the overall robust counter part of the ARO formulation is given as follows 

\begin{equation}
         \min_{\ait, \bit}  \Omega 
\end{equation}
with epigraph constraint  
\begin{align}
\begin{split}
& \sum_{t=1}^{T} \ait + \min_{s^1_1, s^1_2} \sum_{j=1}^{2} c_j s_j \leq \Omega ,\\
& s^1_1 - s^1_2 \geq \bit, \quad \forall t=2,\dots,T ,\\
& s^1_1 - s^1_2 \geq 0 ,\\
& s^1_1, s^1_2 \geq 0 .
\end{split}
\end{align}
With loss coverage constraint 
 \begin{align}
     \begin{split}
 & -\sum_{t=1}^{T} \ait + \sum_{j} c_j s^2_j \leq -\delta , \\ 
& s^2_1 - s^2_2 \geq 1 - \bit, \quad \forall t=2,\dots,T , \\ 
& s^2_1 - s^2_2 \geq 1, \\ 
& s^2_1, s^2_2 \geq 0.
     \end{split}
 \end{align}
With positivity constraint 
 \begin{align}
     \begin{split}
 & \ait + \sum_{j} c_j s^{3,t}_j \leq 0 ,\\ 
& s^{3,t}_1 - s^{3,t}_2 \geq - \beta_{i, t} \\ 
& s^{3,t}_1 - s^{3,t}_2 \geq 0 ,\\ 
& s^{3,t}_1, s^{3,t}_2 \geq 0. 
     \end{split}
 \end{align}
With slowly varying  constraint 
\begin{align}
    \begin{split}
        &      || \alpha_{i,t} - \alpha_{i, t-1} || \le \gamma_3,
    \quad  t = 2, \dots ,T , \\ 
&     || \beta_{i,t} - \beta_{i, t-1} || \le \gamma_4,
    \quad  t = 2, \dots ,T . \\ 
    \end{split}
\end{align}

\section{Case Study for US National Flood Insurance}
In this section, we demonstrate the application of our framework through a case study on flood insurance in the United States. We explain the model parameter estimation and details on the training of machine learning risk predictions. 

\subsection{Data}
We used two redacted datasets from the National Flood Insurance Program (NFIP) on claims and policies respectively. Both data sets are created and maintained by Federal Emergency Management Agency (FEMA). The claims transaction data provide details on NFIP claims transactions across from all states in the United States \cite{nfip_claims}. This dataset consist of 2570089 lines of claim transactions ranging, dated from 1970. 
Due to limited data availability in the early years, we included data from 1975 to 2022. \\ 

For this study, we aggregate data into state level on an annual basis, and have used the following features: date ('dateOfLoss'), state, claim amount ('amountPaidOnBuildingClaim'). Additionally, data from 'MP' (Northern Mariana Islands), 'AS' (American Samoa), 'GU' (Guam), and 'DC' (District of Columbia) have been omitted due to their limited data records. As a result, our cleaned dataset encompasses information from 52 jurisdictions over 48 years, including 50 U.S. states, alongside two territories recognized as island states: the U.S. Virgin Islands and Puerto Rico, enhancing the geographical breadth of our study. \\ 

The policy premium data contains 228664 lines of data
containing policies from 2009 to 2022 \cite{nfip_policy}. We use the policy data as a benchmark to compare model performance between last ten years of testing period from 2013 to 2022. Similar to the claims data, we aggregated data into state level on an annual basis, and have used the following features: state ('propertyState'), date ('policyTeminationDate') and premium ('totalInsurancePremiumOfThePolicy'). \\

In this work, we consider setting insurance premium at state level on an annual basis. We consider setting the premium for the last 10 years, from 2013 to 2022. We trained machine learning models using data from 1975 to 2012, to derive risks in the testing period from 2013 to 2022. We used machine learning derived risks as parameter input for the optimization model, and drive premiums using the RO and ARO frameworks.

\subsection{Optimization Model Parameter Estimation}
Recall from equation \ref{eq:L_CLT}, to compute $L_i^{CLT}$ we need to compute the historical mean and variance for each state. We use the training data from 1975 to 2012 to estimate optimization model parameter to avoid data spoilage in the testing data set. We compute the historical mean and standard deviation for all states on an annual basis. Table \ref{tab:mean_std} in the appendix exhibits the historical mean and variance for the top 10 most costly states. \\ 

In this work, we model the demand sensitivity to insurance premium through a piece-wise linear demand function. 
We estimate the decline rate using historical data from several states. 
We include Figure \ref{fig:demand_la} and \ref{fig:demand_ny} in the appendix showing the scatter plot of number of policy holders in a year against the mean policy premium of that state at that year.  
Different states have different degrees of sensitivity to price, but in general we observe a downward trend of decline in policy holder number as a function of increased price. For the illustrative purpose of this work, we do not specify different sensitivity in different states, but use the same demand damping function across all states. \\

We use the following piece-wise linear demand damping function to model demand damping 

\begin{align} \label{eq:piece_wise_linear}
    f(p) = 
\begin{cases}
    1, & \text{if } p \leq P_0, \\
    1 - m \cdot (p - P_0), & \text{if } p > P_0 \text{ and } f(p) \geq c_{\text{min}}, \\
    c_{\text{min}}, & \text{otherwise}.
\end{cases}
\end{align}

where $P_{0}$ is the minimum premium at which demand damping starts to occur, and $c_{min}$ is the minimum fraction of demand. In this work, we choose $P_0$ to be $\frac{1}{10}\cdot P_{hist}^{max}$, a fraction of maximum historical premium ever charged. We choose $c_{min}$ to be $0.2$ representing at least 20\% of the total demand is preserved regardless of price. We experimented with different demand damping rate: $m = 1/P_{hist}^{max}$. We include Figure \ref{fig:demand_curve} in the appendix to illustrate several choices of the demand damping curve. 

\subsection{Machine Learning Model}

The goal of the machine learning model is to obtain the risk measure $q_{i,k}$, denoting the risk of major flooding event occurrence in the next $k$ years at state $i$, in the machine learning risk uncertainty \ref{eq:ro_ml_u_set}. 
Recall that we denote a major flooding event as the total state-level loss surpassing a certain threshold in the next 1-K years, and obtain the probabilistic prediction result obtained from the binary machine learning task. \\

Specifically, we train machine learning models to predict future annual losses exceeding specific thresholds, based on both current and past losses, across three different threshold levels and time frames. 
We train binary classification models to predict the target, with 1 indicating for a particular state at a particular year, the state will suffer an annual loss passing through the threshold $\Theta$ in the next 1 to K years. We have experimented with three threshold values, corresponding to 90th, 95th and 99th percentile annual claim amount values across all states over all training data, corresponding to USD 18,558,788, USD 50,688,672 and USD 321,903,271. In addition, we have experimented with three K values, corresponding to 3, 5, 10 years respectively. 
A detailed explanation of the data processing, feature construction, and model training methodology is provided in the appendix.\\

Table \ref{tab:ml_consolidated} records out-of-sample prediction results using testing data, corresponding to data between 2012 to 2022. 
We treat data in the testing period on a rolling basis, and we drop the years where we do not have target data, i.e., in year 2019, we predict for K=3 but not for K = 5 or 10. 
We remark that accuracy is generally higher for longer forecasting horizons. 
This is intuitive for the following reason:  given a longer forecasting window, there is a higher probability of a flood occurring for a flood-prone region, leading to a more balanced data. 
Finally, we use the probabilistic prediction results for each state at each testing year $q_{i,k}$ as input to construct uncertainty sets for the robust optimization model as given by equation \ref{eq:L_ML}.

\begin{table}[ht]
\centering
\begin{tabular}{c|cc|cc|cc}
\hline
\multirow{2}{*}{Scores} & \multicolumn{2}{c|}{3 Years} & \multicolumn{2}{c|}{5 Years} & \multicolumn{2}{c}{10 Years} \\ 
 & logreg & xgb & logreg & xgb & logreg & xgb \\ \hline 
\multicolumn{7}{c}{90\% Threshold} \\ \hline
auc & 0.743 & 0.776 & 0.763 & 0.808 & 0.767 & 0.912 \\
f1 & 0.630 & 0.665 & 0.593 & 0.735 & 0.621 & 0.817 \\
accu & 0.693 & 0.675 & 0.632 & 0.736 & 0.623 & 0.830 \\
accu\_bl & 0.625 & 0.706 & 0.605 & 0.753 & 0.641 & 0.809 \\
precision & 0.407 & 0.457 & 0.520 & 0.625 & 0.658 & 0.777 \\
recall & 0.437 & 0.793 & 0.362 & 0.901 & 0.500 & 0.967 \\ \hline
\multicolumn{7}{c}{95\% Threshold} \\ \hline
auc & 0.818 & 0.871 & 0.818 & 0.897 & 0.794 & 0.921 \\
f1 & 0.684 & 0.673 & 0.700 & 0.744 & 0.764 & 0.763 \\
accu & 0.856 & 0.781 & 0.808 & 0.792 & 0.830 & 0.774 \\
accu\_bl & 0.665 & 0.742 & 0.699 & 0.829 & 0.736 & 0.820 \\
precision & 0.300 & 0.306 & 0.368 & 0.460 & 0.595 & 0.560 \\
recall & 0.391 & 0.688 & 0.516 & 0.891 & 0.500 & 0.938 \\ \hline
\multicolumn{7}{c}{99\% Threshold} \\ \hline
auc & 0.920 & 0.945 & 0.922 & 0.956 & 0.952 & 0.992 \\
f1 & 0.704 & 0.687 & 0.737 & 0.771 & 0.835 & 0.906 \\
accu & 0.910 & 0.950 & 0.909 & 0.943 & 0.925 & 0.962 \\
accu\_bl & 0.809 & 0.646 & 0.786 & 0.731 & 0.958 & 0.979 \\
precision & 0.253 & 0.215 & 0.313 & 0.380 & 0.556 & 0.714 \\
recall & 0.696 & 0.304 & 0.640 & 0.480 & 1.000 & 1.000 \\ \hline
\end{tabular}
\caption{Consolidated out-of-sample accuracy for percentile threshold predictions at 90\%, 95\%, and 99\%, across different time horizons using logistic regression (logreg) and XGBoost (xgb) models.}
\label{tab:ml_consolidated}
\end{table}

\section{Results}
We implement three types of robust optimization models, with linearly decreasing demand damping and machine learning risk forecasts. In addition, we implement the adaptive robust optimization model. 
We compare results against two baseline policy premiums. 
\begin{itemize}
    \item Historical premiums charged for each state, which is referred to as `hist'. 
    \item Cumulative moving average loss up to that year, which we refer to as `CMA'. We compute the cumulative moving average loss as follows 
    \begin{equation}
        p^{CMA}_{i,t} = \frac{1}{t}\sum_{t'=0}^t \l_{i, t'}.
    \end{equation}
\end{itemize}
An overall summary of models implemented can be found in Table \ref{tab:model_summary} below. \\ 

\begin{table}[ht]
\centering
\begin{tabular}{l|p{10cm}}
\hline
\textbf{Model Name} & \textbf{Description} \\
\hline
Hist & Uses the historical level of insurance premium collected by the NFIP program, relying on past data to set future premiums without adjustment for future uncertainties. \\
\hline
CMA & Employs a cumulative moving average to compute insurance premiums using historical losses, aiming for stability through past data analysis. \\
\hline
RO1 & Implements robust optimization with a CLT (Central Limit Theorem) uncertainty set, focusing on mitigating risk within a specific range of uncertainty based on statistical theory. \\
\hline
RO2 & Extends robust optimization by incorporating both CLT and ML (Machine Learning) uncertainty sets, aiming for a comprehensive approach to risk management by leveraging advanced analytics. \\
\hline
ARO & Adaptive robust optimization that dynamically adjusts insurance premiums in response to changing conditions and uncertainties, optimizing strategies over time to minimize risk. \\
\hline
\end{tabular}
\caption{Summary of models for insurance premium calculation.}
\label{tab:model_summary}
\end{table}

 We evaluate performance during the testing period, the last ten years of available data in the NFIP data set between 2013 to 2022. For the rest of the section, we choose $\gamma_1$ to be 50000, and $\delta$ to be 10000. We also undertake a sensitivity analysis to assess the impact of our model's parameter selections. Overall, we evaluate our model performance using the following two criteria:
\begin{itemize}
    \item Effectiveness: to evaluate the effectiveness of models to cover losses at a state level. We count the number of insolvent states, i.e., the cumulative premium collected over the testing period does not cover cumulative claim losses. 
    \item Efficiency: to ensure models are charging reasonable levels of premiums to cover losses, we evaluate the overall surplus (or deficit) level as well as the absolute deviation from actual losses to evaluate models' capabilities in realistically assessing risks. 
\end{itemize}


\subsection{Sensitivity Analysis} \label{Sensitivity}
We examine the model sensitivity to parameter choices in our model. Specifically, we examine the effect of the value of $\gamma_2$ and the demand damping rate $m$. Recall that $\gamma_2$ is the parameter controlling how conservative the CLT uncertainty set given by equation \ref{eq:L_CLT}, and the demand damping rate controls how fast demand declines in response to increase in price given by equation \ref{eq:piece_wise_linear}. To examine the overall performance of the premium, we compute the cumulative surplus $S$ across all states over the testing period as follows 
\begin{equation}
      S(\gamma_2) =  \sum_{i=1}^{N} \sum_{t=1}^{T1} \pit - \sum_{i=1}^{N}\sum_{t=1}^{T1} \lit^{act}.
\end{equation}

\begin{figure}[ht]
    \centering
\includegraphics[width=0.75\linewidth]{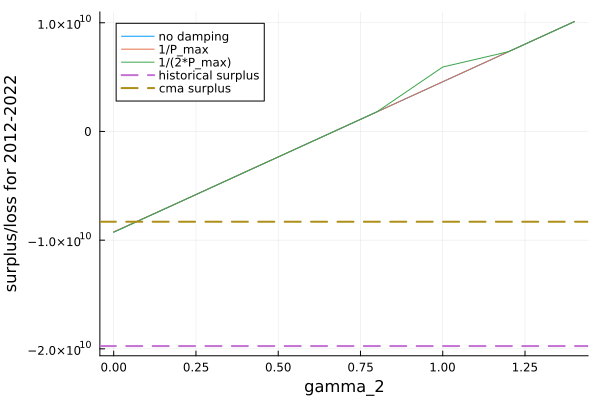}
    \caption{Surplus (or loss) computed during 2012 to 2022 across all states when vary different levels of $\gamma_2$. Two dotted lines demonstrate the level of surplus from two baseline models: historical surplus calculated from the actual premiums collected, cma surplus is computed using the cumulative moving average. }
    \label{fig:vary_gamma2}
\end{figure}

Figure \ref{fig:vary_gamma2} shows the level of surplus as a function of different $\gamma_2$, with different demand damping rate. The two dotted line shows the constant surplus computed by two baselines: using the actual premiums collected during this period and the CMA rule. We observe that both baselines incur a loss over the testing period, with historical premiums resulting in about 20 billion loss, and CMA rule resulting in 8 million loss. \\

We let $\gamma_2$ to take values between 0 and 1.5 at a stepsize of 0.1, and resolve the optimization model at each $\gamma_2$ value, and compute the sum of surplus across all states across testing years. $\gamma_2 = 0$ corresponding to convex optimization without uncertainty, and $\gamma_2 = 1.5$ with the maximum degree of uncertainty. As we increase the value of $\gamma_2$, the size of the uncertainty set increases, and the model becomes more conservative resulting in surplus as expected. We remark that the surplus breaks even when $\gamma_2$ takes value between 0.6 and 0.7. And we observe a smooth increase in surplus as $\gamma_2$ increases. \\ 

In addition, we experiment with three demand damping rates: no damping, $m_1=1/P_{hist}^{max}$, $m_2=1/(2*P_{hist}^{max})$, with $m_2$ damps twice as fast as $m_1$. Similar as above, we experiment with varying $\gamma_2$ corresponding to the different demand damping rates. We observe that the choice of demand damping is less significant compared with the variation of $\gamma_2$. 
  
\subsection{Effectiveness} 


Evaluating the effectiveness of insurance schemes in covering losses is crucial, particularly through the lens of insolvency rates. Therefore, we consider the insolvency status at a state level, where the cumulative premium fails to cover the cumulative loss, as an indicator of a scheme's financial resilience and risk management efficiency. Our analysis spans the 2013 - 2022 testing period, corresponding the the last ten years of data, and focuses on the impact of $\gamma_2$, a parameter which significantly influences model outcomes. As detailed in Table \ref{tab:insolvency_rates}, we explore the effects of varying $\gamma_2$ values, from 0 to 2 in increments of 0.2, on the number of insolvent states induced by each model.\\ 

Note that the performance of the Historical and CMA schemes is invariant to changes in $\gamma_2$. Specifically, the CMA approach results in 36 insolvent states, outperforming the Historical scheme, which results in insolvency across all states. This outcome underscores the limitations of relying solely on past data for setting future premiums. In contrast, both the RO and ARO schemes exhibit increased conservatism—and consequently, fewer insolvent states—with higher $\gamma_2$ values. This is because $\gamma_2$ directly relates to the size of the CLT uncertainty set, as specified by equation \ref{eq:L_CLT}, where larger $\gamma_2$ values necessitate higher premium values to mitigate risk.\\ 

This mechanism within the RO framework affords decision-makers the flexibility to adjust the conservatism of their models based on risk tolerance and financial strategy. Among the evaluated schemes, the RO2 model, which integrates both CLT and ML uncertainty sets, consistently achieves the lowest number of insolvent states for a given $\gamma_2$ level. This performance is followed closely by the ARO scheme and the RO1 model. Such findings highlight the nuanced balance between risk management and financial sustainability, illustrating the advanced capabilities of RO2 and ARO in navigating the complexities of insurance premium optimization under uncertainty.

\begin{table}[h]
\centering
\begin{tabular}{c|c|c|c|c|c}
\hline
\textbf{$\gamma_2$} & \textbf{ARO} & \textbf{RO1} & \textbf{RO2} & \textbf{CMA} & \textbf{Hist} \\
\hline
0.0 & 36.0 & 36.0 & 32.0 & 36.0 & 52.0 \\
0.2 & 30.0 & 32.0 & 28.0 & 36.0 & 52.0 \\
0.4 & 26.0 & 26.0 & 23.0 & 36.0 & 52.0 \\
0.6 & 25.0 & 26.0 & 23.0 & 36.0 & 52.0 \\
0.8 & 22.0 & 23.0 & 20.0 & 36.0 & 52.0 \\
1.0 & 21.0 & 22.0 & 20.0 & 36.0 & 52.0 \\
1.2 & 19.0 & 21.0 & 19.0 & 36.0 & 52.0 \\
1.4 & 16.0 & 18.0 & 16.0 & 36.0 & 52.0 \\
1.6 & 15.0 & 15.0 & 14.0 & 36.0 & 52.0 \\
1.8 & 13.0 & 13.0 & 12.0 & 36.0 & 52.0 \\
2.0 & 12.0 & 12.0 & 11.0 & 36.0 & 52.0 \\
\hline
\end{tabular}
\caption{The number of insolvent states during the testing period based on different methods. In total we test over 52 states, which are the 50 US states with additionally Puerto Rico and the US virgin islands.}
\label{tab:insolvency_rates}
\end{table}

\subsection{Efficiency}
In addition, we consider the efficiency of the models by considering the surplus and deficit level, as well as the absolute deviation from the actual loss incurred. This is to ensure models are not over-charging states, and pricing premiums correctly align with risks. \\

Table \ref{tab:combined_insolvency_surplus_deficit} exhibits the overall surplus (or deficit) level and the absolute deviation level across all states over the testing period. We compute the overall absolute deviation across all states over the testing period from the actual loss, which is an indicator of how well models are able to trace the risks. We compute AD as a function of $\gamma_2$ as follows
\begin{equation}
    AD(\gamma_2) =  \sum_{i=1}^{N} \sum_{t=1}^{T1} |\pit - \lit^{act}|.
\end{equation}

\begin{figure}[h]
    \centering
    \includegraphics[width=0.85\textwidth]{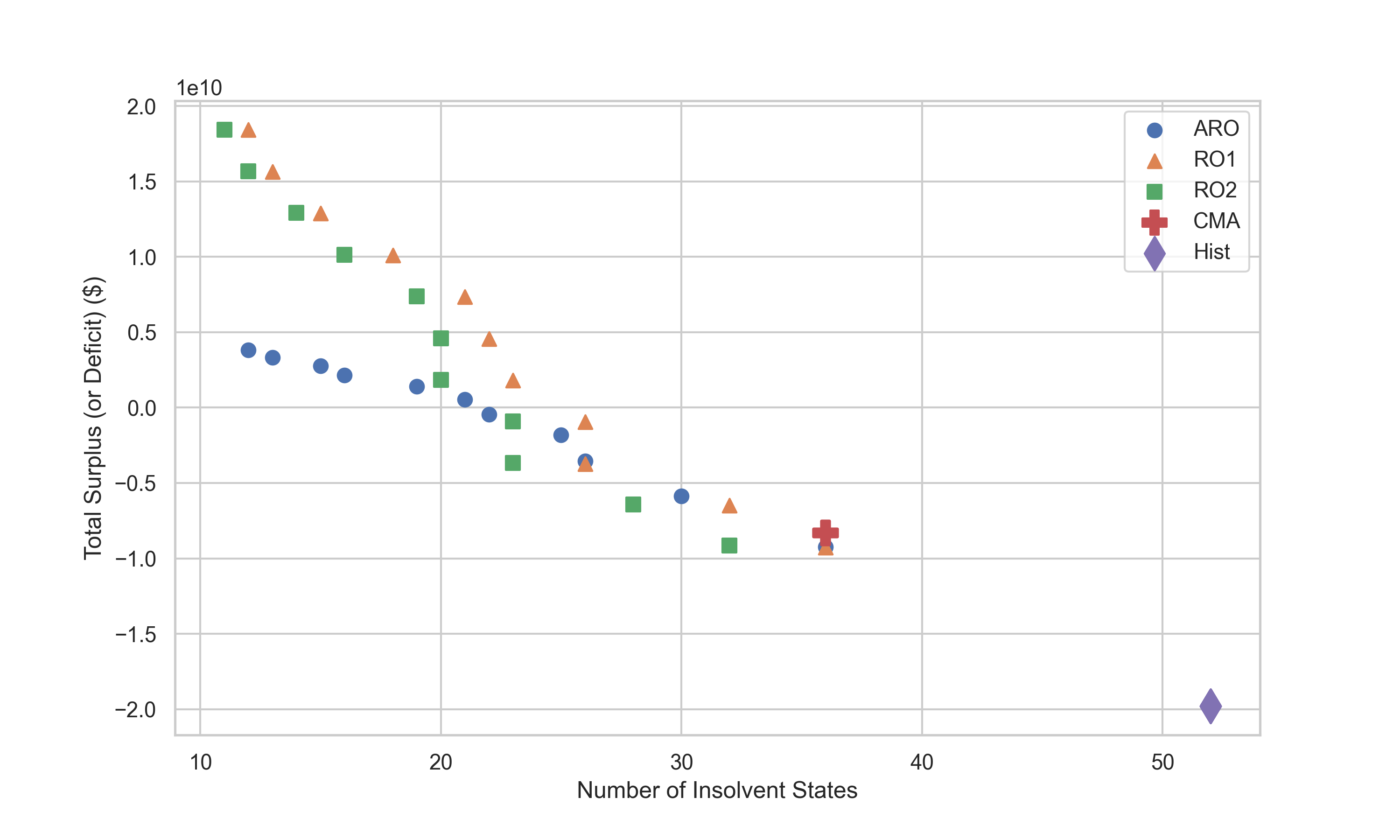}
    \caption{
    The scatter plot visualizes the efficient frontier, showing how different values of $\gamma_2$ affect the number of insolvent states (x-axis) and the total surplus (or deficit) (y-axis) computed as the total premium charged minus actual loss over the testing period. Note that CMA and Hist are plotted as static points because their values do not change with varying $\gamma_2$ values. }
    \label{fig:efficient_frontier}
\end{figure}

First, we observe that Historical premiums significantly undercharges over the testing period, resulting in 19 billion losses, and CMA rule results in 8 million losses. This suggests historical levels are insufficient to cover future losses. Similar as before, we observe as $\gamma_2$ increases, RO and ARO schemes increase the level of conservatism and reaches more surplus. In particular, comparing the two RO models, RO2 scheme is more conservative than RO1, albeit not significantly, as expected. Because RO2 contains the additioanl ML uncertainty set. With same level of $\gamma_2$, we observe that ARO achieves the surplus the slower than RO schemes. Especially when $\gamma_2$ reaches the level of 1.4, ARO scheme increases premiums much slower than RO schemes, which increases at a constant rate with increasing $\gamma_2$. Figure \ref{fig:efficient_frontier} illustrates the efficient frontier illustrating the trade-off between the number of insolvent states versus the surplus (or deficit) achieved during the testing period.\\

The superiority of the ARO scheme over RO scheme is reinforced when looking at the absolute deviation metric. We observe that the level of error using the ARO scheme is more stable than both RO schemes. The stability of ARO scheme makes it desirable for policy-makers, as the scheme offers great robustness against insolvency with the least premium charged when chosen a high level of $\gamma_2$. In conclusion, the results underscore the ARO and RO models' superior performance in managing uncertainty and adapting to changing risk profiles, highlighting their potential for enhancing catastrophe insurance premium pricing strategies. We include Figure \ref{fig:efficient_frontier2} in the appendix illustrating the trade-off between the number of insolvent states versus the absolute deviation during the testing period. \\

\begin{table}[ht]
\centering
\resizebox{\textwidth}{!}{%
\begin{tabular}{c|cc|cc|cc|cc|cc}
\hline
\textbf{$\gamma_2$} & \multicolumn{2}{c|}{\textbf{ARO}} & \multicolumn{2}{c|}{\textbf{RO1}} & \multicolumn{2}{c|}{\textbf{RO2}} & \multicolumn{2}{c|}{\textbf{CMA}} & \multicolumn{2}{c}{\textbf{Hist}} \\
 & S/D & A/D & S/D & A/D & S/D & A/D & S/D & A/D & S/D & A/D \\
\hline
0.0 & $-9.23e9$ & $2.26e10$ & $-9.27e9$ & $2.25e10$ & $-9.16e9$ & $2.26e10$ & $-8.31e9$ & $2.30e10$ & $-1.98e10$ & $1.98e10$ \\
0.2 & $-5.88e9$ & $2.47e10$ & $-6.50e9$ & $2.42e10$ & $-6.41e9$ & $2.42e10$ & $-8.31e9$ & $2.30e10$ & $-1.98e10$ & $1.98e10$ \\
0.4 & $-3.56e9$ & $2.65e10$ & $-3.73e9$ & $2.59e10$ & $-3.67e9$ & $2.59e10$ & $-8.31e9$ & $2.30e10$ & $-1.98e10$ & $1.98e10$ \\
0.6 & $-1.80e9$ & $2.78e10$ & $-0.96e9$ & $2.77e10$ & $-0.92e9$ & $2.77e10$ & $-8.31e9$ & $2.30e10$ & $-1.98e10$ & $1.98e10$ \\
0.8 & $-0.46e9$ & $2.88e10$ & $1.80e9$  & $2.95e10$ & $1.84e9$  & $2.95e10$ & $-8.31e9$ & $2.30e10$ & $-1.98e10$ & $1.98e10$ \\
1.0 & $0.55e9$  & $2.96e10$ & $4.57e9$  & $3.13e10$ & $4.60e9$  & $3.13e10$ & $-8.31e9$ & $2.30e10$ & $-1.98e10$ & $1.98e10$ \\
1.2 & $1.41e9$  & $3.03e10$ & $7.34e9$  & $3.32e10$ & $7.37e9$  & $3.32e10$ & $-8.31e9$ & $2.30e10$ & $-1.98e10$ & $1.98e10$ \\
1.4 & $2.16e9$  & $3.09e10$ & $10.11e9$ & $3.53e10$ & $10.14e9$ & $3.53e10$ & $-8.31e9$ & $2.30e10$ & $-1.98e10$ & $1.98e10$ \\
1.6 & $2.78e9$  & $3.13e10$ & $12.88e9$ & $3.74e10$ & $12.90e9$ & $3.74e10$ & $-8.31e9$ & $2.30e10$ & $-1.98e10$ & $1.98e10$ \\
1.8 & $3.32e9$  & $3.18e10$ & $15.64e9$ & $3.96e10$ & $15.67e9$ & $3.96e10$ & $-8.31e9$ & $2.30e10$ & $-1.98e10$ & $1.98e10$ \\
2.0 & $3.82e9$  & $3.21e10$ & $18.41e9$ & $4.17e10$ & $18.43e9$ & $4.17e10$ & $-8.31e9$ & $2.30e10$ & $-1.98e10$ & $1.98e10$ \\
\hline
\end{tabular}
} 
\caption{Condensed table showing surplus/deficit (S/D) and absolute deviation (A/D) across different $\gamma_2$ levels for ARO, RO1, RO2, CMA, and Hist.}
\label{tab:combined_insolvency_surplus_deficit}
\end{table}

\section{Conclusion}

In conclusion, we present an Adaptive Robust Optimization (ARO) framework to catastrophe insurance premium pricing. We first present a nominal linear optimization formulation to introduce the problem of setting insurance prices for rare catastrophe events, and present a robust optimization formulation with two distinct uncertainty sets. The Central Limit Theorem (CLT) uncertainty set protects against deviations from historical losses, and the Machine Learning (ML) uncertainty set to incorporate predicted risks. We derive the robust counter part by solving the inner problem in closed form, and present a convex optimization re-formulation. In addition, we extend the RO framework to an Adaptive Robust Optimization (ARO) model, using linear decision rules and enabling premium adjustments based on actual loss experiences. \\

We applied the framework to the US flood insurance and evaluate our performance against two baseline benchmarks: the US National Flood Insurnace Program (NFIP) premiums and CMA premiums. We employed historical data from 1975 to 2012 to construct the uncertainty sets and train machine learning models, and we used the last ten years of available NFIP data to evaluate RO insurance scheme against the two benchmarks: the actual historical premium policies and premiums derived from cumulative moving average (CMA) rules. We demonstrate the superiority of an ARO approach in two metrics: effectiveness and efficiency. First, optimization-based models are able to effectively cover losses, achieving a smooth transition from high number of insolvent states to a low number, thus granting policy providers the discretion to determine the desired insolvency level. Second, optimization-based models are capable of efficiently charge premiums, resulting in a smooth transition from deficit to surplus balance depending on model parameter value. In particular, we recommend policy makers to use an ARO model, with conservative parameter values, to achieve superior performance in both effectiveness and efficiency, resulting in achieving simultaneously low insolvent rate and relatively low premiums charged.

\newpage 
\bibliographystyle{plain}
\bibliography{insurance}

\newpage 
\appendix 
\section{Optimization Model Parameter Estimation}
In this section, we expand on how to estimate parameters for the optimization model. 

\subsection{Computing $L_i^{CLT}$}
Recall from equation \ref{eq:L_CLT}, to compute $L_i^{CLT}$ we need to compute the historical mean and variance for each state. Table \ref{tab:mean_std} in the appendix exhibiting historical mean and standard deviation for the top 10 most costly states, computed on the annual basis. 
\begin{table}[h] 
\centering 
\begin{tabular}{lcccc}
\hline
State & Max & Mean & Std & Median \\
\hline
LA & 3,763,390 & 45,084 & 58,467 & 20,583 \\
TX & 8,973,270 & 44,046 & 65,420 & 19,801 \\
NJ & 4,022,518 & 32,727 & 58,000 & 13,042 \\
NY & 9,467,720 & 35,725 & 74,510 & 13,306 \\
FL & 9,100,033 & 25,351 & 63,809 & 8,052 \\
MS & 10,000,000 & 46,603 & 94,924 & 14,823 \\
NC & 1,294,678 & 21,725 & 41,345 & 7,992 \\
PA & 1,889,793 & 19,040 & 41,943 & 6,939 \\
AL & 4,900,000 & 28,059 & 95,474 & 8,466 \\
SC & 1,764,000 & 25,574 & 43,752 & 10,088 \\
\hline
\end{tabular}
\caption{Statistics of top 10 most costly states in the US, including maximum annual claim loss, mean, standard deviation and median.}
\label{tab:mean_std}
\end{table}

\subsection{Demand Damping}
In this work, we model the demand sensitivity to insurance premium through a piece-wise linear demand function. We estimate the decline rate using historical data from several states. Figures \ref{fig:demand_la} and \ref{fig:demand_ny} below shows the scatter plot of number of policy holders in a year against the mean policy premium of that state at that year.  Different states have different degrees of sensitivity to price, but in general we observe a downward trend of decline in policy holder number as a function of increased price. For the illustrative purpose of this work, we do not specify different sensitivity in different states, but use the same demand damping function across all states. \\
\begin{figure}[h]
  \centering
  \begin{minipage}[b]{0.45\textwidth}
    \centering
    \includegraphics[width=\textwidth]{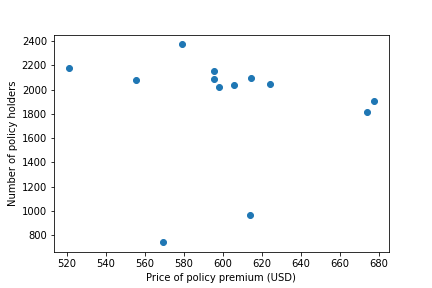}
    \caption{(a) Demand damping estimation for Louisiana state (LA).}
    \label{fig:demand_la}
  \end{minipage}
  \hfill
  \begin{minipage}[b]{0.45\textwidth}
    \centering
    \includegraphics[width=\textwidth]{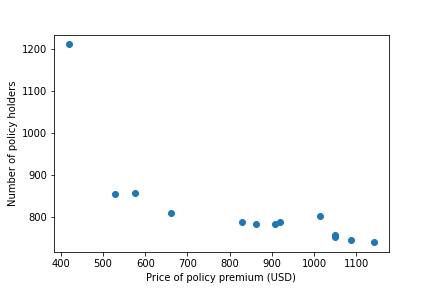}
    \caption{(b) Demand damping estimation for New York state (NY).}
    \label{fig:demand_ny}
  \end{minipage}
  \end{figure}

\begin{figure}
    \centering
    \includegraphics[width=0.6\linewidth]{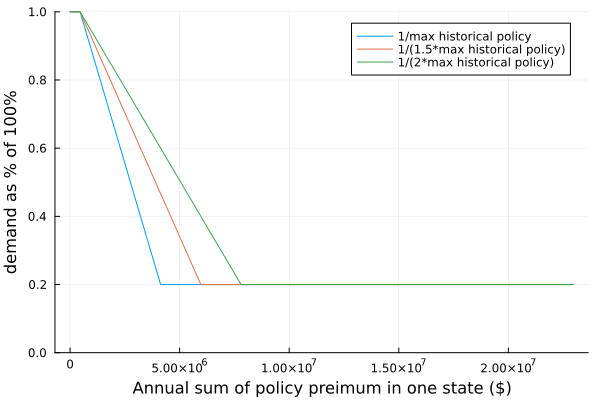}
    \caption{Different piece-wise linear demand damping curves corresponding to different rate of decline.}
    \label{fig:demand_curve}
\end{figure}

\newpage 
\section{Detailed Training Protocol of Machine Learning Models}
\subsection{Data Processing}

Our dataset spans from 1975 to 2022. We exclude any claims lacking a specified claim amount to ensure the integrity of our analysis. Additionally, data from 'MP' (Northern Mariana Islands), 'AS' (American Samoa), 'GU' (Guam), and 'DC' (District of Columbia) have been omitted due to their limited data records. As a result, our cleaned dataset encompasses information from 52 jurisdictions over 48 years. This includes all 50 U.S. states, alongside two territories recognized as island states: the U.S. Virgin Islands and Puerto Rico, enhancing the geographical breadth of our study.\\ 

To construct the machine learning model, first, we aggregate data to state and annual level. The index of the data is two levels: state and year. 
For missing data, we performed linear interpolation within each state using previous and later years. 
Then for each state at a particular year, we construct the following features: state (categorical), current year annual loss, past 1-5 years annual loss. Current and past year losses are numerical features, where as state (categorical) feature is treated with one-hot encoding.  \\ 

We train binary classification models to predict the target, with 1 indicating for a particular state at a particular year, the state will suffer an annual loss passing through the threshold $\Theta$ in the next 1 to K years. We have experimented with three threshold values, corresponding to 90th, 95th and 99th percentile annual claim amount values across all states over all training data, corresponding to USD 18,558,788, USD 50,688,672 and USD 321,903,271. In addition, we have experimented with three K values, corresponding to 3, 5, 10 years respectively.\\

\subsection{Training and Testing Protocols}
We split the data set chronologically into training period from 1975 to 2011, and testing period from 2012 to 2022. We experiment with two standard machine learning models, logistic regression and XGBoost. We employed 3-fold cross validation to search for the best parameters for each type of models. The search space for hyperparameter tuning can be found in Table \ref{tab:hyperparameters} below. \\

\begin{table}[h]
\centering
\begin{tabular}{c|c|c}
\hline
\textbf{Model} & \textbf{Hyperparameters} & \textbf{Values} \\
\hline \hline 
XGBoost &
number of estimators & 100, 150 \\
& maximum tree depth & 4, 6\\
& learning rate & 0.1, 0.3 \\
\hline
Logistic Regression &
C & 0, 0.2, 0.4, 0.6, 0.8, 1\\
& penalty & L1, L2 \\
\hline
\end{tabular}
\caption{Hyperparameters searched for our models.}
\label{tab:hyperparameters}
\end{table}

\subsection{Detailed Prediction Results}
\begin{figure}
    \centering
    \includegraphics[width=0.6\textwidth]{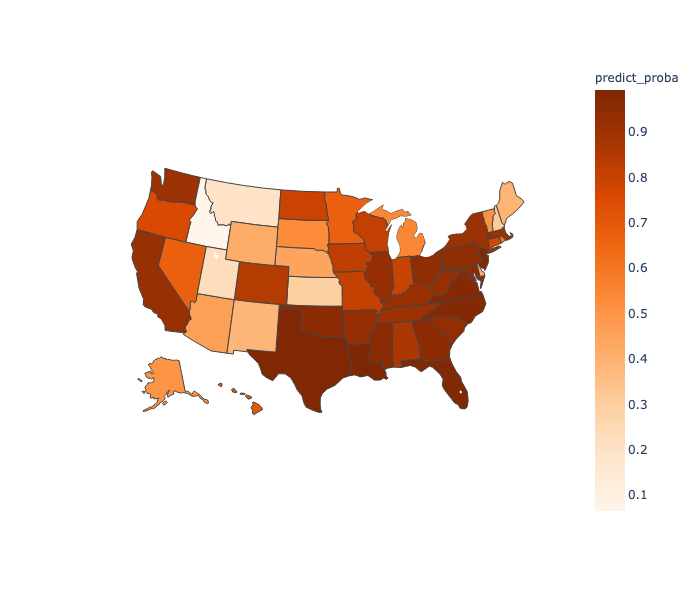}
    \caption{Map illustrating machine learning (ML) predicted risks for all states surpassing 90th-percentile flooding risk within a 5-year time frame from 2016. Regions are color-coded, with darker shades indicating a higher probability of predicted risks.}
    \label{fig:ml_predictions}
\end{figure}

Figure \ref{fig:ml_predictions} illustrates machine learning predicted risks for all states surpassing 90th percentile flooding risk within the next 5-year time frame on 2016. For each state, we produce one prediction for each year over the testing period, for each threshold, and for 3-year, 5-year, 10-year time frames. 
Table \ref{tab:ml_consolidated} record out-of-sample prediction results using testing data, corresponding to data between 2012 to 2022. 
We treat data in the testing period on a rolling basis, and we drop the years where we do not have target data, i.e., in year 2019, we predict for K=3 but not for K = 5 or 10. 
We remark that accuracy is generally higher for longer forecasting horizons. This is likely due to the following reasons: first, longer forecasting horizon lead to higher probability of flood, which leads to more balanced data; second, we have less testing samples. 
We use the probabilistic prediction results for each state at each testing year $q_{i,k}$ as input to construct uncertainty sets for the robust optimization model as given by equation \ref{eq:L_ML}. 

\subsection{Computational resources}
The source codes for this study, implemented in Julia 1.7 and Python 3.9, are publicly accessible at [repository link]. The convex optimization problems were solved using Ipopt and Gurobi solvers, with machine learning models trained on a local machine with 4 Intel CPUs on a Macbook Pro personal computer. Comprehensive documentation detailing the methodology and specific parameters can be found in the repository's code comments.

\section{Additional Plot of the Efficient Frontier}
\begin{figure}[h]
    \centering
    \includegraphics[width=0.8\textwidth]{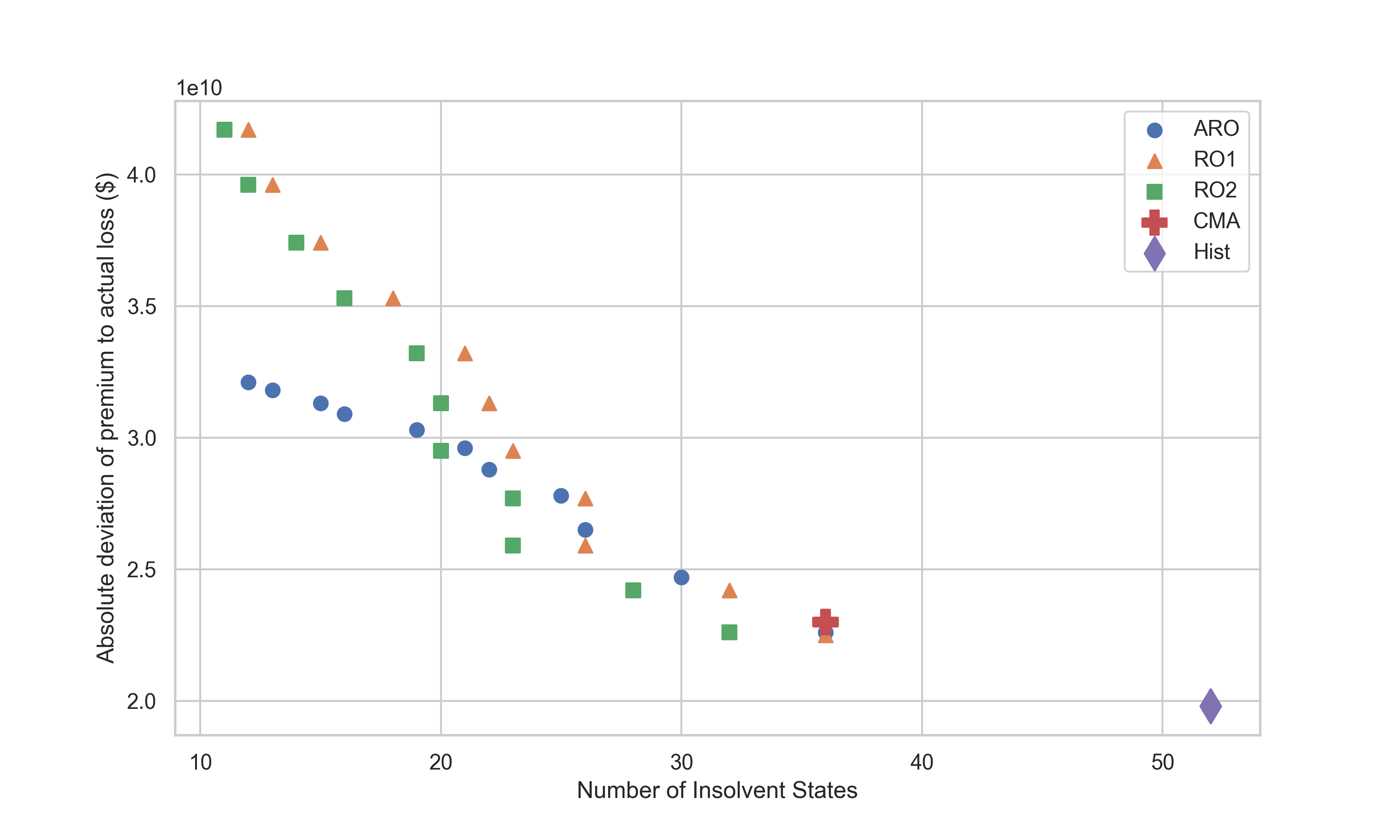}
    \caption{
    The scatter plot visualizes the efficient frontier, showing how different values of $\gamma_2$ affect the number of insolvent states (x-axis) and the total surplus (or deficit) (y-axis) computed as the total premium charged minus actual loss over the testing period. Note that CMA and Hist are plotted as static points because their values do not change with varying $\gamma_2$ values. The plot demonstrates that ARO achieves higher efficiency with lower absolute deviation at high $\gamma_2$ values, highlighting its better performance under these conditions.}
    \label{fig:efficient_frontier2}
\end{figure}

Figure \ref{fig:efficient_frontier2} illustrates the trade-off between the number of insolvent states versus the absolute deviation during the testing period. 

\end{document}